\documentclass[12pt]{article}
\usepackage[reqno]{amsmath}
\usepackage{amssymb, theorem, enumerate}
\usepackage{tikz}
\usepackage{url}

\expandafter\ifx\csname pdfoptionalwaysusepdfpagebox\endcsname\relax\else
\fi

\theoremstyle{change}
{\theorembodyfont{\slshape}
\newtheorem{theorem}{Theorem.}[section]
\newtheorem{lemma}[theorem]{Lemma.}
}
\theorembodyfont{\rmfamily}

\newcommand\eps{\epsilon}

\newcommand\cref[1]{Corollary~\ref{cor:#1}}

\newcommand\ff[1]{{\langle x \rangle_#1}}

\def\proof{\noindent{{\sl Proof. }}}
\def\sqr#1#2{{\vbox{\hrule height.#2pt
    \hbox{\vrule width.#2pt height#1pt \kern#1pt
        \vrule width.#2pt}\hrule height.#2pt}}}
\def\eqed{\sqr53}
\def\qed{%
    \ifmmode\eqno\eqed
    \else\nobreak\ \hfill\eqed\medbreak\fi}

\title{Planar triangulations with real chromatic roots arbitrarily close to four}
\author{Gordon Royle\\
	School of Computer Science \& Software Engineering\\
	University of Western Australia\\
	 35 Stirling Highway\\
	 Nedlands, WA 6009, Australia.}

\begin{document}
\maketitle

\begin{abstract}
We exhibit infinite families of planar graphs with real chromatic roots 
arbitrarily close to 4, thus resolving a long-standing conjecture in the affirmative.
\end{abstract}

\section{Introduction}

The chromatic polynomial $P(G,x)$ of a graph $G$ is the function that counts the number of proper $x$-colourings of a graph when $x$ is a positive integer. It was introduced by Birkhoff \cite{MR1502436} in 1912, in the hope that a quantitative study of the {\em numbers} of colourings of planar graphs could be used to show that $P(G,4) > 0$, thus resolving the 4-colour conjecture. Although defined on positive integers, it is well-known that $P(G,x)$ is a polynomial in $x$, and therefore we can view it as a function on the real or complex numbers.  Although they were unable to prove the 4-colour conjecture with chromatic polynomials, Birkhoff and Lewis \cite{MR0018401} showed that for $G$ planar, $P(G,x) > 0$ for all $x \in [5,\infty)$, and they conjectured that $P(G,x) > 0$ for all $x \in [4,\infty)$, thus providing the first results and problems in the theory of {\em real} chromatic roots.  

With the 4-colour conjecture still unsolved, the study of the {\em complex} chromatic roots of planar graphs was initiated in the 1960s (Hall, Siry \& Vanderslice \cite{MR0179110}, Berman \& Tutte \cite{MR0238734}) when the availability of computers made empirical study of complex chromatic roots feasible. A decade later, at around the same time as the 4-colour conjecture became a theorem with a decidedly {\em combinatorial} proof, Beraha \& Kahane \cite{MR539072} found families of planar graphs with complex chromatic roots arbitrarily close to 4, thus providing strong evidence that there could never be a complex-analytic proof of the 4-colour theorem. 

Despite this, the study of real and complex chromatic roots of not-necessarily-planar graphs has developed into a substantial theory, significantly enhanced by the relationship between the chromatic polynomial and the partition function of the $q$-state Potts model of statistical mechanics. Jackson \cite{MR2005532} provides an excellent and comprehensive survey of the main results and questions in the theory of chromatic roots, while an accessible introduction to the connections between the chromatic polynomial and the Potts model can be found in Sokal \cite{alan99}.

One of the unresolved problems described by Jackson \cite{MR2005532} was a conjecture of Beraha which (among other things) asserted the existence of planar graphs with {\em real} chromatic roots arbitrarily close to 4. At the time of Jackson's survey the largest known real chromatic root of a planar graph was $3.8267\ldots$ obtained from a 21-vertex graph found by Woodall. In this paper we resolve this conjecture in the affirmative, describing a number of related families of planar triangulations with real chromatic roots arbitrarily close to 4. 

\section{Double-ended lattice graphs}

\begin{figure}
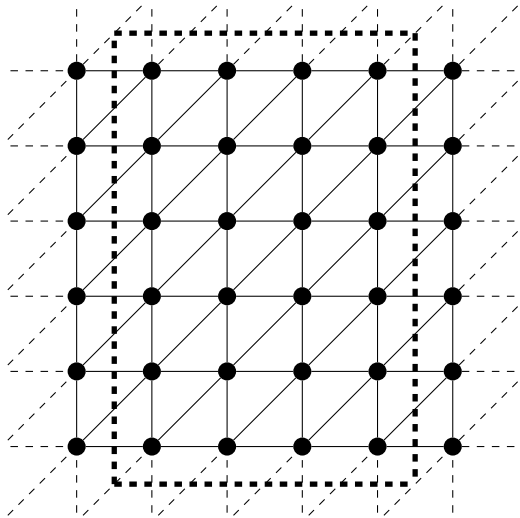

\begin{center}
\begin{pgfpicture}{0cm}{0cm}{7cm}{7cm}
\color{white}
\pgfnodecircle{v00}[fill]{\pgfxy(0,0)}{0.12cm}
\pgfnodecircle{v01}[fill]{\pgfxy(0,1)}{0.12cm}
\pgfnodecircle{v02}[fill]{\pgfxy(0,2)}{0.12cm}
\pgfnodecircle{v03}[fill]{\pgfxy(0,3)}{0.12cm}
\pgfnodecircle{v04}[fill]{\pgfxy(0,4)}{0.12cm}
\pgfnodecircle{v05}[fill]{\pgfxy(0,5)}{0.12cm}
\pgfnodecircle{v06}[fill]{\pgfxy(0,6)}{0.12cm}
\pgfnodecircle{v07}[fill]{\pgfxy(0,7)}{0.12cm}
\pgfnodecircle{v10}[fill]{\pgfxy(1,0)}{0.12cm}
\pgfnodecircle{v17}[fill]{\pgfxy(1,7)}{0.12cm}
\pgfnodecircle{v20}[fill]{\pgfxy(2,0)}{0.12cm}
\pgfnodecircle{v27}[fill]{\pgfxy(2,7)}{0.12cm}
\pgfnodecircle{v30}[fill]{\pgfxy(3,0)}{0.12cm}
\pgfnodecircle{v37}[fill]{\pgfxy(3,7)}{0.12cm}
\pgfnodecircle{v40}[fill]{\pgfxy(4,0)}{0.12cm}
\pgfnodecircle{v47}[fill]{\pgfxy(4,7)}{0.12cm}
\pgfnodecircle{v50}[fill]{\pgfxy(5,0)}{0.12cm}
\pgfnodecircle{v57}[fill]{\pgfxy(5,7)}{0.12cm}
\pgfnodecircle{v60}[fill]{\pgfxy(6,0)}{0.12cm}
\pgfnodecircle{v67}[fill]{\pgfxy(6,7)}{0.12cm}
\pgfnodecircle{v70}[fill]{\pgfxy(7,0)}{0.12cm}
\pgfnodecircle{v71}[fill]{\pgfxy(7,1)}{0.12cm}
\pgfnodecircle{v72}[fill]{\pgfxy(7,2)}{0.12cm}
\pgfnodecircle{v73}[fill]{\pgfxy(7,3)}{0.12cm}
\pgfnodecircle{v74}[fill]{\pgfxy(7,4)}{0.12cm}
\pgfnodecircle{v75}[fill]{\pgfxy(7,5)}{0.12cm}
\pgfnodecircle{v76}[fill]{\pgfxy(7,6)}{0.12cm}
\pgfnodecircle{v77}[fill]{\pgfxy(7,7)}{0.12cm}
\color{black}
\pgfnodecircle{v11}[fill]{\pgfxy(1,1)}{0.12cm}
\pgfnodecircle{v12}[fill]{\pgfxy(1,2)}{0.12cm}
\pgfnodecircle{v13}[fill]{\pgfxy(1,3)}{0.12cm}
\pgfnodecircle{v14}[fill]{\pgfxy(1,4)}{0.12cm}
\pgfnodecircle{v15}[fill]{\pgfxy(1,5)}{0.12cm}
\pgfnodecircle{v16}[fill]{\pgfxy(1,6)}{0.12cm}
\pgfnodecircle{v21}[fill]{\pgfxy(2,1)}{0.12cm}
\pgfnodecircle{v22}[fill]{\pgfxy(2,2)}{0.12cm}
\pgfnodecircle{v23}[fill]{\pgfxy(2,3)}{0.12cm}
\pgfnodecircle{v24}[fill]{\pgfxy(2,4)}{0.12cm}
\pgfnodecircle{v25}[fill]{\pgfxy(2,5)}{0.12cm}
\pgfnodecircle{v26}[fill]{\pgfxy(2,6)}{0.12cm}
\pgfnodecircle{v31}[fill]{\pgfxy(3,1)}{0.12cm}
\pgfnodecircle{v32}[fill]{\pgfxy(3,2)}{0.12cm}
\pgfnodecircle{v33}[fill]{\pgfxy(3,3)}{0.12cm}
\pgfnodecircle{v34}[fill]{\pgfxy(3,4)}{0.12cm}
\pgfnodecircle{v35}[fill]{\pgfxy(3,5)}{0.12cm}
\pgfnodecircle{v36}[fill]{\pgfxy(3,6)}{0.12cm}
\pgfnodecircle{v41}[fill]{\pgfxy(4,1)}{0.12cm}
\pgfnodecircle{v42}[fill]{\pgfxy(4,2)}{0.12cm}
\pgfnodecircle{v43}[fill]{\pgfxy(4,3)}{0.12cm}
\pgfnodecircle{v44}[fill]{\pgfxy(4,4)}{0.12cm}
\pgfnodecircle{v45}[fill]{\pgfxy(4,5)}{0.12cm}
\pgfnodecircle{v46}[fill]{\pgfxy(4,6)}{0.12cm}
\pgfnodecircle{v51}[fill]{\pgfxy(5,1)}{0.12cm}
\pgfnodecircle{v52}[fill]{\pgfxy(5,2)}{0.12cm}
\pgfnodecircle{v53}[fill]{\pgfxy(5,3)}{0.12cm}
\pgfnodecircle{v54}[fill]{\pgfxy(5,4)}{0.12cm}
\pgfnodecircle{v55}[fill]{\pgfxy(5,5)}{0.12cm}
\pgfnodecircle{v56}[fill]{\pgfxy(5,6)}{0.12cm}
\pgfnodecircle{v61}[fill]{\pgfxy(6,1)}{0.12cm}
\pgfnodecircle{v62}[fill]{\pgfxy(6,2)}{0.12cm}
\pgfnodecircle{v63}[fill]{\pgfxy(6,3)}{0.12cm}
\pgfnodecircle{v64}[fill]{\pgfxy(6,4)}{0.12cm}
\pgfnodecircle{v65}[fill]{\pgfxy(6,5)}{0.12cm}
\pgfnodecircle{v66}[fill]{\pgfxy(6,6)}{0.12cm}
\pgfnodeconnline{v11}{v21}
\pgfnodeconnline{v11}{v12}
\pgfnodeconnline{v11}{v22}
\pgfnodeconnline{v12}{v22}
\pgfnodeconnline{v12}{v13}
\pgfnodeconnline{v12}{v23}
\pgfnodeconnline{v13}{v23}
\pgfnodeconnline{v13}{v14}
\pgfnodeconnline{v13}{v24}
\pgfnodeconnline{v14}{v24}
\pgfnodeconnline{v14}{v15}
\pgfnodeconnline{v14}{v25}
\pgfnodeconnline{v15}{v25}
\pgfnodeconnline{v15}{v16}
\pgfnodeconnline{v15}{v26}
\pgfnodeconnline{v16}{v26}
\pgfnodeconnline{v21}{v31}
\pgfnodeconnline{v21}{v22}
\pgfnodeconnline{v21}{v32}
\pgfnodeconnline{v22}{v32}
\pgfnodeconnline{v22}{v23}
\pgfnodeconnline{v22}{v33}
\pgfnodeconnline{v23}{v33}
\pgfnodeconnline{v23}{v24}
\pgfnodeconnline{v23}{v34}
\pgfnodeconnline{v24}{v34}
\pgfnodeconnline{v24}{v25}
\pgfnodeconnline{v24}{v35}
\pgfnodeconnline{v25}{v35}
\pgfnodeconnline{v25}{v26}
\pgfnodeconnline{v25}{v36}
\pgfnodeconnline{v26}{v36}
\pgfnodeconnline{v31}{v41}
\pgfnodeconnline{v31}{v32}
\pgfnodeconnline{v31}{v42}
\pgfnodeconnline{v32}{v42}
\pgfnodeconnline{v32}{v33}
\pgfnodeconnline{v32}{v43}
\pgfnodeconnline{v33}{v43}
\pgfnodeconnline{v33}{v34}
\pgfnodeconnline{v33}{v44}
\pgfnodeconnline{v34}{v44}
\pgfnodeconnline{v34}{v35}
\pgfnodeconnline{v34}{v45}
\pgfnodeconnline{v35}{v45}
\pgfnodeconnline{v35}{v36}
\pgfnodeconnline{v35}{v46}
\pgfnodeconnline{v36}{v46}
\pgfnodeconnline{v41}{v51}
\pgfnodeconnline{v41}{v42}
\pgfnodeconnline{v41}{v52}
\pgfnodeconnline{v42}{v52}
\pgfnodeconnline{v42}{v43}
\pgfnodeconnline{v42}{v53}
\pgfnodeconnline{v43}{v53}
\pgfnodeconnline{v43}{v44}
\pgfnodeconnline{v43}{v54}
\pgfnodeconnline{v44}{v54}
\pgfnodeconnline{v44}{v45}
\pgfnodeconnline{v44}{v55}
\pgfnodeconnline{v45}{v55}
\pgfnodeconnline{v45}{v46}
\pgfnodeconnline{v45}{v56}
\pgfnodeconnline{v46}{v56}
\pgfnodeconnline{v51}{v61}
\pgfnodeconnline{v51}{v52}
\pgfnodeconnline{v51}{v62}
\pgfnodeconnline{v52}{v62}
\pgfnodeconnline{v52}{v53}
\pgfnodeconnline{v52}{v63}
\pgfnodeconnline{v53}{v63}
\pgfnodeconnline{v53}{v54}
\pgfnodeconnline{v53}{v64}
\pgfnodeconnline{v54}{v64}
\pgfnodeconnline{v54}{v55}
\pgfnodeconnline{v54}{v65}
\pgfnodeconnline{v55}{v65}
\pgfnodeconnline{v55}{v56}
\pgfnodeconnline{v55}{v66}
\pgfnodeconnline{v56}{v66}
\pgfnodeconnline{v61}{v62}
\pgfnodeconnline{v62}{v63}
\pgfnodeconnline{v63}{v64}
\pgfnodeconnline{v64}{v65}
\pgfnodeconnline{v65}{v66}
\pgfsetdash{{0.1cm}{0.1cm}}{0cm}

\pgfnodeconnline{v00}{v11}
\pgfnodeconnline{v01}{v11}

\pgfnodeconnline{v01}{v12}
\pgfnodeconnline{v02}{v12}

\pgfnodeconnline{v02}{v13}
\pgfnodeconnline{v03}{v13}

\pgfnodeconnline{v03}{v14}
\pgfnodeconnline{v04}{v14}

\pgfnodeconnline{v04}{v15}
\pgfnodeconnline{v05}{v15}

\pgfnodeconnline{v05}{v16}
\pgfnodeconnline{v06}{v16}

\pgfnodeconnline{v10}{v11}
\pgfnodeconnline{v10}{v21}
\pgfnodeconnline{v16}{v17}
\pgfnodeconnline{v16}{v27}

\pgfnodeconnline{v20}{v21}
\pgfnodeconnline{v20}{v31}
\pgfnodeconnline{v26}{v27}
\pgfnodeconnline{v26}{v37}

\pgfnodeconnline{v30}{v31}
\pgfnodeconnline{v30}{v41}
\pgfnodeconnline{v36}{v37}
\pgfnodeconnline{v36}{v47}

\pgfnodeconnline{v40}{v41}
\pgfnodeconnline{v40}{v51}
\pgfnodeconnline{v46}{v47}
\pgfnodeconnline{v46}{v57}

\pgfnodeconnline{v50}{v51}
\pgfnodeconnline{v50}{v61}
\pgfnodeconnline{v56}{v57}
\pgfnodeconnline{v56}{v67}

\pgfnodeconnline{v60}{v61}

\pgfnodeconnline{v61}{v71}
\pgfnodeconnline{v61}{v72}
\pgfnodeconnline{v62}{v72}
\pgfnodeconnline{v62}{v73}
\pgfnodeconnline{v63}{v73}
\pgfnodeconnline{v63}{v74}
\pgfnodeconnline{v64}{v74}
\pgfnodeconnline{v64}{v75}
\pgfnodeconnline{v65}{v75}
\pgfnodeconnline{v65}{v76}
\pgfnodeconnline{v66}{v76}
\pgfnodeconnline{v66}{v67}
\pgfnodeconnline{v66}{v77}

\pgfsetlinewidth{2pt}
\pgfline{\pgfxy(1.5,0.5)}{\pgfxy(5.5,0.5)}
\pgfline{\pgfxy(1.5,6.5)}{\pgfxy(5.5,6.5)}
\pgfline{\pgfxy(1.5,0.5)}{\pgfxy(1.5,6.5)}
\pgfline{\pgfxy(5.5,6.5)}{\pgfxy(5.5,0.5)}
\end{pgfpicture}
\end{center}
\caption{A $4 \times 6$ strip of the infinite triangular lattice}
\label{fig:lat}
\end{figure}

Some of the most important families of graphs from a statistical mechanics point of view are rectangular subsets of infinite two-dimensional lattices, particularly the square and triangular lattices (see Figure~\ref{fig:lat} for the triangular lattice).

Physicists are usually interested in the {\em limiting behaviour} of the complex chromatic roots as the strip becomes infinitely wide and tall. In practice, this is studied by fixing one of the parameters, say the width, and examining the curves where the chromatic roots accumulate as the height tends to infinity. Repeating this calculation for a number of different fixed widths then gives insight into how the entire curves move as the width increases. In addition, they consider the cases of {\em periodic boundary conditions} where the left edge of the strip is viewed as being adjacent to the right edge as though the strip were wrapped around a cylinder, and {\em free boundary conditions} otherwise. These different boundary conditions are denoted by using $P$ and $F$ respectively as subscripts, so for example $m_P \times n_F$ denotes an $m \times n$ rectangular strip (of whichever type of lattice is under consideration) with periodic transverse (i.e. left-right) boundary conditions and free longitudinal (i.e. top-bottom) boundary conditions.

For families of graphs of this form with one parameter fixed, the limiting behaviour of the chromatic roots can be explained by a theorem of Beraha, Kahane \& Weiss \cite{MR520560} which shows that the accumulation points of the roots form one or more curves and isolated points in the complex plane. The location and nature of these curves has been heavily studied by a number of researchers, in particular the prolific group of Shrock, Chang \& Tsai who, with various collaborators, have studied both the full Tutte polynomial  (a.k.a ``the partition function of the $q$-state Potts model'') and the chromatic polynomial for various different lattices with a number of different boundary conditions. Shrock \cite{MR1821978} gives a good introduction to their earlier work on chromatic polynomials, while the recent work \cite{chsh05} contains many other references. For the {\em triangular} lattice, which is the focus of this paper, the full Tutte polynomial for strips of this lattice has been computed and studied for widths up to 5 by Chang, Jacobsen, Salas \& Shrock \cite{MR2035629}, and the chromatic polynomial for widths up to 10 by Jacobsen, Salas \& Sokal \cite{MR2000228}.

The graphs found by Beraha \& Kahane \cite{MR539072} with complex chromatic roots arbitrarily close to 4 are $4_P \times n_F$ strips of the triangular lattice with an additional vertex joined to each end. Ro{\v{c}}ek, Shrock \& Tsai \cite{MR1625635, roshts} considered the effect of adding a number of different end graphs to either one end or both ends of a strip, and it is these {\em double-ended} lattice graphs that we can use to find real chromatic roots arbitrarily close to 4.

In particular we show that if we attach {\em suitable pairs} of planar end-graphs to the top and bottom of the cylindrical triangular lattice $4_P \times n_F$, then the resulting families of graphs have {\em real} chromatic roots arbitrarily close to 4.

\section{Chromatic Polynomials}

Let $A$ and $B$ be two graphs each with a distinguished 4-cycle, say $a_1a_2a_3a_4$ and $b_1b_2b_3b_4$ respectively, and suppose that $G$ is the graph obtained by {\em gluing together} $A$ and $B$, i.e. identifying $a_i$ with $b_i$ for $1 \leq i \leq 4$. Any proper colouring of $G$ determines proper colourings of $A$ and $B$ such that the distinguished 4-cycle is coloured identically in both graphs. Conversely, if we have enough information about the colourings of $A$ and $B$ {\em and} the colourings of the distinguished 4-cycle in both $A$ and $B$, then we can determine the chromatic polynomial of $G$.

Given a graph $A$ with distinguished 4-cycle $a_1a_2a_3a_4$, we can assign one of four types to each proper colouring $\varphi$ of $A$ according to the colour partition it induces on the 4-cycle, as follows:

\begin{description}
\item[Type 1] $\varphi(a_1) = \varphi(a_3)$ and $\varphi(a_2) = \varphi(a_4)$,
\item[Type 2] $\varphi(a_1) = \varphi(a_3)$ and $\varphi(a_2) \not= \varphi(a_4)$,
\item[Type 3] $\varphi(a_1) \not= \varphi(a_3)$ and $\varphi(a_2) = \varphi(a_4)$,
\item[Type 4] $\varphi(a_1) \not= \varphi(a_3)$ and $\varphi(a_2) \not= \varphi(a_4)$.
\end{description}

The {\em number} of colourings of each type is just the chromatic polynomial of an auxiliary graph obtained from $A$ by identifying two vertices if they are required to have the same colour, and adding an edge between them if they are required to have different colours. Let $P_1(A,x)$, $P_2(A,x)$, $P_3(A,x)$ and $P_4(A,x)$ denote the four chromatic polynomials corresponding to the four classes above, and note that
$$
P(A,x) = P_1(A,x)+P_2(A,x)+P_3(A,x)+P_4(A,x).
$$
We will find it convenient to express this as a vector in ${\mathbb Z}[x]^4$ and so we define the {\em partitioned chromatic polynomial} $Q(A,x)$ of graph with a distinguished 4-cycle to be
$$
Q(A,x) = 
\left(
\begin{array}{c}
P_1(A,x)\\
P_2(A,x)\\
P_3(A,x)\\
P_4(A,x)
\end{array}
\right).
$$

\begin{lemma}\label{lem:glue}
Let $A$ and $B$ be two graphs with distinguished 4-cycles and corresponding partitioned chromatic polynomials $Q(A,x)$ and  $Q(B,x)$. Then the chromatic polynomial of the graph $G$ obtained by gluing together $A$ and $B$ is the sole
entry of the $1 \times 1$ matrix
$$
Q(A)^T  D Q(B)
$$
where
$$ 
D = 
\left(
\begin{array}{cccc}
1/ \ff 2 & 0 & 0 & 0\\
0 & 1/ \ff 3 & 0 & 0\\
0 & 0 & 1 / \ff 3 & 0\\
0 & 0 & 0 & 1 / \ff 4
\end{array}
\right),
$$
and $\ff k$ denotes the $k$'th {\em falling factorial} $x(x-1)\cdots (x-k+1)$.
\end{lemma}
 
\proof 
For a positive integer $x$, the product $P_i(A,x) P_i(B,x)$ counts the number of pairs $(\rho,\varphi)$ such that $\rho$ is a proper $x$-colouring of type $i$ of $A$, and $\varphi$ is a proper $x$-colouring of type $i$ of $B$. The two colourings can be combined to form a colouring of $G$ if and only if they both use the same colours on the distinguished 4-cycle, and all proper $x$-colourings of $G$ arise in this fashion. If the colourings of type $i$ use $s$ colours on the distinguished 4-cycle, then the fraction of these pairs of colourings using the same colours on the distinguished 4-cycle is exactly $1/\ff s$. Summing the four terms of this form, and expressing the result in matrix terms gives the stated result.
\qed

\begin{figure}
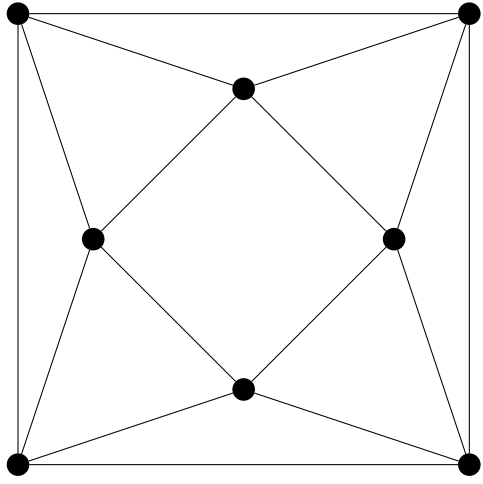

\begin{center}
\begin{pgfpicture}{0cm}{0cm}{6cm}{6cm}
\pgfnodecircle{v0}[fill]{\pgfxy(0,0)}{0.15cm}
\pgfnodecircle{v1}[fill]{\pgfxy(6,0)}{0.15cm}
\pgfnodecircle{v2}[fill]{\pgfxy(6,6)}{0.15cm}
\pgfnodecircle{v3}[fill]{\pgfxy(0,6)}{0.15cm}

\pgfnodecircle{v8}[fill]{\pgfxy(3,1)}{0.15cm}
\pgfnodecircle{v9}[fill]{\pgfxy(5,3)}{0.15cm}
\pgfnodecircle{v10}[fill]{\pgfxy(3,5)}{0.15cm}
\pgfnodecircle{v11}[fill]{\pgfxy(1,3)}{0.15cm}

\pgfnodeconnline{v0}{v1}
\pgfnodeconnline{v1}{v2}
\pgfnodeconnline{v2}{v3}
\pgfnodeconnline{v3}{v0}

\pgfnodeconnline{v0}{v8}
\pgfnodeconnline{v8}{v1}
\pgfnodeconnline{v1}{v9}
\pgfnodeconnline{v9}{v2}
\pgfnodeconnline{v2}{v10}
\pgfnodeconnline{v10}{v3}
\pgfnodeconnline{v3}{v11}
\pgfnodeconnline{v11}{v0}

\pgfnodeconnline{v9}{v10}
\pgfnodeconnline{v10}{v11}
\pgfnodeconnline{v8}{v9}
\pgfnodeconnline{v11}{v8}

\end{pgfpicture}
\end{center}
\caption{$L$ is one layer of the cylindrical triangular lattice}
\label{fig:gadget}
\end{figure}

Now consider the ``gadget'' $L$ shown in Figure~\ref{fig:gadget}; this is actually a $4_P \times 2_F$ strip of the triangular lattice, or in other words a single ``layer'' of the cylindrical triangular lattice of width 4. Let $M$ denote the matrix whose rows and columns are indexed by the 4 types of colour partition of a 4-cycle and where  the entry $M_{ij}$ is the polynomial counting the number of $x$-colourings of $L$ that are of type $i$ on the outer 4-cycle, and 
type $j$ on the inner 4-cycle.  Every partition of $L$ into $s$ independent sets contributes a term of $\ff s$ to one of the entries of $M$, and a simple computation yields 
{\small
$$
M = 
{ 
\left(
\begin{array}{cccc}
\ff 4  & \ff 5  & \ff 5  & \ff 6\\
\ff 5  & \ff 4 + 2 \ff 5 + \ff 6 & \ff  4 + 2 \ff 5 + \ff 6 & 4 \ff 5 + 4 \ff 6 + \ff 7 \\
\ff 5 & \ff 4 + 2 \ff 5 + \ff 6 & \ff  4 + 2 \ff 5 + \ff 6 & 4 \ff 5 + 4 \ff 6 + \ff 7  \\
\ff 6 & 4 \ff 5 + 4 \ff 6 + \ff 7 & 4 \ff 5 + 4 \ff 6 + \ff 7 & M_{44}
\end{array}
\right)
}
$$
}
where
$$
M_{44} = 
2 \ff 4 + 16 \ff 5 + 20 \ff 6 + 8 \ff 7 + \ff 8.
$$

Now we consider the effect of adding a layer of the cylindrical triangular lattice to an existing graph. In particular, take a graph $A$ with a distinguished 4-cycle, and identify this 4-cycle with the inner 4-cycle of $L$, forming a new graph $A'$ whose distinguished 4-cycle is the outer 4-cycle of $L$. In order to repeat this operation with $A'$, we need some way to ``transfer'' the colouring information from the inner 4-cycle to the outer 4-cycle.

\begin{lemma}\label{lem:gadget}
Let $A$ be a graph with a distinguished 4-cycle, and suppose that $A'$ is obtained by identifying the inner 4-cycle of $L$ with the distinguished 4-cycle of $A$, and declaring the outer 4-cycle of $L$ to be the distinguished 4-cycle of $A'$. Then the partitioned chromatic polynomial $Q(A',x)$ is given by
$$
Q(A') = M D Q(A).
$$
\end{lemma}

\proof 
This is a counting argument analogous to that of Lemma~\ref{lem:glue}, where the proper colourings of $A'$ are counted according to the type of colour partition on the inner 4-cycle of $L$, and then allocated to the result according to the type of colour partition induced on the outer 4-cycle of $L$.
 \qed

The matrix $MD$ is the {\em transfer matrix} for the cylindrical triangular lattice, because it transfers the colouring information from one layer of the lattice to the next. In one form or another transfer-matrix techniques are the main tool for studying chromatic and Tutte polynomials of lattice graphs --- more details on two approaches to using transfer matrices can be found in Chang \& Shrock \cite{chsh05} and Salas \& Sokal \cite{MR1853428}.

\begin{theorem}\label{thm:main}
Let $A$ and $B$ be graphs with distinguished 4-cycles, and let $X_{A,B}(n)$ be the graph obtained from the triangular lattice strip $4_P \times n_F$ by gluing $A$ to the 4-cycle at the top of the strip and $B$ to the 4-cycle at the bottom. Then the chromatic polynomial of $X_{A,B}(n)$ is the sole entry of the $1 \times 1$ matrix
$$
Q(A)^T D (MD)^{n-1} Q(B).
$$
\end{theorem}

\proof
This follows directly from Lemma~\ref{lem:glue} and Lemma~\ref{lem:gadget}.
\qed

\begin{figure}
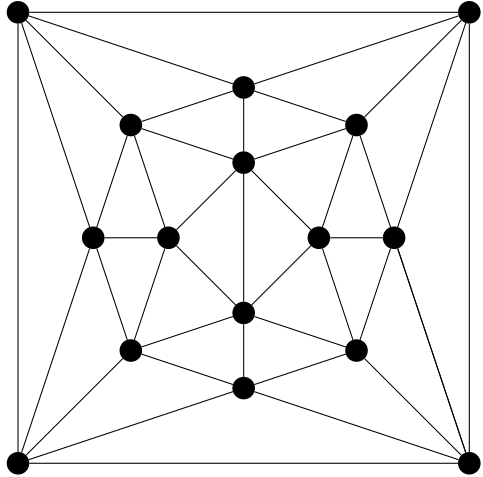

\begin{center}
\begin{pgfpicture}{0cm}{0cm}{6cm}{6cm}
\pgfnodecircle{v0}[fill]{\pgfxy(0,0)}{0.15cm}
\pgfnodecircle{v1}[fill]{\pgfxy(6,0)}{0.15cm}
\pgfnodecircle{v2}[fill]{\pgfxy(6,6)}{0.15cm}
\pgfnodecircle{v3}[fill]{\pgfxy(0,6)}{0.15cm}
\pgfnodecircle{v4}[fill]{\pgfxy(1.5,1.5)}{0.15cm}
\pgfnodecircle{v5}[fill]{\pgfxy(4.5,1.5)}{0.15cm}
\pgfnodecircle{v6}[fill]{\pgfxy(4.5,4.5)}{0.15cm}
\pgfnodecircle{v7}[fill]{\pgfxy(1.5,4.5)}{0.15cm}
\pgfnodecircle{v8}[fill]{\pgfxy(3,1)}{0.15cm}
\pgfnodecircle{v9}[fill]{\pgfxy(5,3)}{0.15cm}
\pgfnodecircle{v10}[fill]{\pgfxy(3,5)}{0.15cm}
\pgfnodecircle{v11}[fill]{\pgfxy(1,3)}{0.15cm}
\pgfnodecircle{v12}[fill]{\pgfxy(3,2)}{0.15cm}
\pgfnodecircle{v13}[fill]{\pgfxy(4,3)}{0.15cm}
\pgfnodecircle{v14}[fill]{\pgfxy(3,4)}{0.15cm}
\pgfnodecircle{v15}[fill]{\pgfxy(2,3)}{0.15cm}
\pgfnodeconnline{v0}{v1}
\pgfnodeconnline{v1}{v2}
\pgfnodeconnline{v2}{v3}
\pgfnodeconnline{v3}{v0}
\pgfnodeconnline{v0}{v8}
\pgfnodeconnline{v8}{v1}
\pgfnodeconnline{v1}{v9}
\pgfnodeconnline{v9}{v2}
\pgfnodeconnline{v2}{v10}
\pgfnodeconnline{v10}{v3}
\pgfnodeconnline{v3}{v11}
\pgfnodeconnline{v11}{v0}
\pgfnodeconnline{v0}{v4}
\pgfnodeconnline{v1}{v5}
\pgfnodeconnline{v2}{v6}
\pgfnodeconnline{v3}{v7}
\pgfnodeconnline{v1}{v9}
\pgfnodeconnline{v4}{v8}
\pgfnodeconnline{v8}{v5}
\pgfnodeconnline{v5}{v9}
\pgfnodeconnline{v9}{v6}
\pgfnodeconnline{v6}{v10}
\pgfnodeconnline{v10}{v7}
\pgfnodeconnline{v7}{v11}
\pgfnodeconnline{v11}{v4}
\pgfnodeconnline{v8}{v12}
\pgfnodeconnline{v9}{v13}
\pgfnodeconnline{v10}{v14}
\pgfnodeconnline{v11}{v15}
\pgfnodeconnline{v4}{v12}
\pgfnodeconnline{v12}{v5}
\pgfnodeconnline{v5}{v13}
\pgfnodeconnline{v13}{v6}
\pgfnodeconnline{v6}{v14}
\pgfnodeconnline{v14}{v7}
\pgfnodeconnline{v7}{v15}
\pgfnodeconnline{v15}{v4}
\pgfnodeconnline{v12}{v13}
\pgfnodeconnline{v13}{v14}
\pgfnodeconnline{v14}{v15}
\pgfnodeconnline{v15}{v12}
\pgfnodeconnline{v12}{v14}
\end{pgfpicture}
\end{center}
\caption{$H$ is a triangulation of a square}
\label{fig:base16}
\end{figure}

\begin{figure}
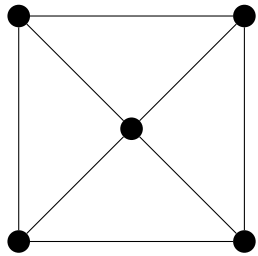

\begin{center}
\begin{pgfpicture}{0cm}{0cm}{3cm}{3cm}
\pgfnodecircle{v0}[fill]{\pgfxy(0,0)}{0.15cm}
\pgfnodecircle{v1}[fill]{\pgfxy(3,0)}{0.15cm}
\pgfnodecircle{v2}[fill]{\pgfxy(3,3)}{0.15cm}
\pgfnodecircle{v3}[fill]{\pgfxy(0,3)}{0.15cm}
\pgfnodecircle{v4}[fill]{\pgfxy(1.5,1.5)}{0.15cm}
\pgfnodeconnline{v0}{v1}
\pgfnodeconnline{v1}{v2}
\pgfnodeconnline{v2}{v3}
\pgfnodeconnline{v3}{v0}
\pgfnodeconnline{v4}{v0}
\pgfnodeconnline{v4}{v1}
\pgfnodeconnline{v4}{v2}
\pgfnodeconnline{v4}{v3}
\end{pgfpicture}
\end{center}
\caption{The 4-wheel $W_4 = C_4 + K_1$}\label{fig:simple}
\end{figure}

Finally, we observe that Woodall's graph with real chromatic root $3.8267\ldots$ is the graph $X_{H,\,W_4}(2)$ where $H$ and $W_4 = C_4 + K_1$ are the graphs shown in Figure~\ref{fig:base16} and Figure~\ref{fig:simple} respectively with the outer face being the distinguished 4-cycle in each case. 

\section{The family $X_{H,\,W_4}(n)$}

In this section we give computational results on the real roots of the family
$X_{H,\,W_4}(n)$, which although not quite the smallest family having real chromatic roots arbitrarily close to 4, is the first family for which this was noticed. It is immediate that
\[
Q(W_4,x) = 
\left(
\begin{array}{c}
\ff3 \\
\ff4 \\
\ff4 \\
\ff5 \\
\end{array}
\right),
\]
and $Q(H)$ is shown in Table~\ref{tab:qb}, where each term is given as the product of a simple multiplier and an explicitly listed polynomial of degree at most 11.

\begin{table}
\begin{center}

\begin{tabular}{|c|rrrr|}
\hline
&\multicolumn{1}{c}{$P_1$}&\multicolumn{1}{c}{$P_2$}&\multicolumn{1}{c}{$P_3$}&\multicolumn{1}{c|}{$P_4$}\\
\hline
Multiplier&\multicolumn{1}{c}{$\ff5$}&\multicolumn{1}{c}{$\ff4$}&\multicolumn{1}{c}{$\ff4$}&\multicolumn{1}{c|}{$\ff4(x-3)$}\\
\hline
$x^{11}$&$$&$1$&$1$&$1$\\
$x^{10}$&$$&$-34$&$-34$&$-34$\\
$x^{9}$&$1$&$538$&$538$&$538$\\
$x^{8}$&$-28$&$-5244$&$-5244$&$-5246 $\\
$x^{7}$&$362 $&$35078 $&$35078 $&$35148 $\\
$x^{6}$&$-2846$&$-169490 $&$-169490 $&$-170548 $\\
$x^{5}$&$15036$&$604806 $&$604806 $&$613920 $\\
$x^{4}$&$-55448$&$-1595807 $&$-1595807 $&$-1645010$\\
$x^{3}$&$142716 $&$3051803 $&$3051803 $&$3222645$\\
$x^{2}$&$-246724 $&$-4024676 $&$-4024676 $&$-4397342$\\
$x^{1}$&$258889 $&$3286881 $&$3286881 $&$3753346$\\
$1$&$-124884$&$-1255163$&$-1255163$&$-1511254 $\\
\hline
\end{tabular}

\end{center}
\caption{Components of the partitioned chromatic polynomial of $H$.}
\label{tab:qb}
\end{table}

Using the expression given in Theorem~\ref{thm:main} and a suitable computer algebra package, it takes only minutes to calculate the  chromatic polynomials for $X_{H,\,W_4}(1)$ to $X_{H,\,W_4}(100)$. The roots of these polynomials were then computed and the largest real roots listed in Table~\ref{tab:reals}.

\begin{table}
\begin{center}
\begin{tabular}{|c|c|}
\hline
$n$&Max. real root for $X_{H,\,W_4}(n)$\\
\hline
1&3.7924699360\\
2&3.8267852044\\
3&3.8483432574\\
4&3.8637744449\\
5&3.8756040984\\
10&3.9100811222\\
20&3.9388450668\\
30&3.9524053758\\
40&3.9605533349\\
50&3.9660736920\\
60&3.9700969536\\
70&3.9731775927\\
80&3.9756223799\\
90&3.9776159809\\
100&3.9792767496\\
\hline
\end{tabular}
\end{center}
\caption{Largest real chromatic roots of $X_{H,\,W_4}(n)$.}
\label{tab:reals}
\end{table}

\begin{figure}
\begin{center}
\begin{tikzpicture}
\draw [step=0.2cm,lightgray,very thin] (0,0) grid (10,14);
\draw [step=2cm,gray,thin,xshift=-1cm,yshift=-1cm] (1,1) grid (11,15);

\foreach \x in {3,5,7,9}
\draw (\x,6.9) -- (\x, 7.1);

\foreach \x in {1,3,5,9,11,13}
\draw (0.9,\x) -- (1.10,\x);

\foreach \x/\xtext in {3/1, 5/2, 7/3, 9/4}
\draw (\x,6.8) node[anchor=north,fill=white] {{\scriptsize $\xtext$}};

\foreach \x/\xtext in {1/-3, 3/-2, 5/-1}
\draw (0.8,\x) node[anchor=east,fill=white] {{\scriptsize $\xtext$}};

\foreach \x/\xtext in {9/1, 11/2, 13/3}
\draw (0.8,\x) node[anchor=east,fill=white] {{\scriptsize $\xtext$}};

\draw (0,7) -- (10,7);
\draw (1,0) -- (1,14);

\draw (3.000,7.000) circle (1pt);
\draw (2.972,1.543) circle (1pt);
\draw (2.972,12.457) circle (1pt);
\draw (3.742,1.498) circle (1pt);
\draw (3.742,12.502) circle (1pt);
\draw (3.747,12.502) circle (1pt);
\draw (3.747,1.498) circle (1pt);
\draw (3.754,1.497) circle (1pt);
\draw (3.754,12.503) circle (1pt);
\draw (3.763,12.503) circle (1pt);
\draw (3.763,1.497) circle (1pt);
\draw (3.776,1.496) circle (1pt);
\draw (3.776,12.504) circle (1pt);
\draw (3.791,12.506) circle (1pt);
\draw (3.791,1.494) circle (1pt);
\draw (3.809,1.493) circle (1pt);
\draw (3.809,12.507) circle (1pt);
\draw (3.830,12.508) circle (1pt);
\draw (3.830,1.492) circle (1pt);
\draw (3.853,1.490) circle (1pt);
\draw (3.853,12.510) circle (1pt);
\draw (3.879,12.511) circle (1pt);
\draw (3.879,1.489) circle (1pt);
\draw (3.908,12.513) circle (1pt);
\draw (3.908,1.487) circle (1pt);
\draw (3.939,12.514) circle (1pt);
\draw (3.939,1.486) circle (1pt);
\draw (3.972,12.516) circle (1pt);
\draw (3.972,1.484) circle (1pt);
\draw (4.009,1.483) circle (1pt);
\draw (4.009,12.517) circle (1pt);
\draw (4.047,12.518) circle (1pt);
\draw (4.047,1.482) circle (1pt);
\draw (4.089,1.481) circle (1pt);
\draw (4.089,12.519) circle (1pt);
\draw (4.132,12.520) circle (1pt);
\draw (4.132,1.480) circle (1pt);
\draw (4.178,1.480) circle (1pt);
\draw (4.178,12.520) circle (1pt);
\draw (4.226,12.520) circle (1pt);
\draw (4.226,1.480) circle (1pt);
\draw (4.277,1.481) circle (1pt);
\draw (4.277,12.519) circle (1pt);
\draw (4.329,12.518) circle (1pt);
\draw (4.329,1.482) circle (1pt);
\draw (4.384,12.516) circle (1pt);
\draw (4.384,1.484) circle (1pt);
\draw (4.441,12.514) circle (1pt);
\draw (4.441,1.486) circle (1pt);
\draw (4.500,12.511) circle (1pt);
\draw (4.500,1.489) circle (1pt);
\draw (4.561,1.493) circle (1pt);
\draw (4.561,12.507) circle (1pt);
\draw (4.623,12.502) circle (1pt);
\draw (4.623,1.498) circle (1pt);
\draw (4.688,1.504) circle (1pt);
\draw (4.688,12.496) circle (1pt);
\draw (4.754,12.490) circle (1pt);
\draw (4.754,1.510) circle (1pt);
\draw (4.822,1.518) circle (1pt);
\draw (4.822,12.482) circle (1pt);
\draw (4.892,1.527) circle (1pt);
\draw (4.892,12.473) circle (1pt);
\draw (4.963,12.463) circle (1pt);
\draw (4.963,1.537) circle (1pt);
\draw (5.000,7.000) circle (1pt);
\draw (5.035,1.548) circle (1pt);
\draw (5.035,12.452) circle (1pt);
\draw (5.109,12.440) circle (1pt);
\draw (5.109,1.560) circle (1pt);
\draw (5.185,1.574) circle (1pt);
\draw (5.185,12.426) circle (1pt);
\draw (5.261,12.412) circle (1pt);
\draw (5.261,1.588) circle (1pt);
\draw (5.339,12.395) circle (1pt);
\draw (5.339,1.605) circle (1pt);
\draw (5.418,1.623) circle (1pt);
\draw (5.418,12.377) circle (1pt);
\draw (5.498,12.358) circle (1pt);
\draw (5.498,1.642) circle (1pt);
\draw (5.579,12.338) circle (1pt);
\draw (5.579,1.662) circle (1pt);
\draw (5.661,1.685) circle (1pt);
\draw (5.661,12.315) circle (1pt);
\draw (5.743,12.292) circle (1pt);
\draw (5.743,1.708) circle (1pt);
\draw (5.826,1.734) circle (1pt);
\draw (5.826,12.266) circle (1pt);
\draw (5.910,12.240) circle (1pt);
\draw (5.910,1.760) circle (1pt);
\draw (5.994,12.212) circle (1pt);
\draw (5.994,1.788) circle (1pt);
\draw (6.079,1.817) circle (1pt);
\draw (6.079,12.183) circle (1pt);
\draw (6.163,12.154) circle (1pt);
\draw (6.163,1.846) circle (1pt);
\draw (6.236,7.000) circle (1pt);
\draw (6.247,1.874) circle (1pt);
\draw (6.247,12.126) circle (1pt);
\draw (6.325,1.904) circle (1pt);
\draw (6.325,12.096) circle (1pt);
\draw (6.397,1.940) circle (1pt);
\draw (6.397,12.060) circle (1pt);
\draw (6.473,1.985) circle (1pt);
\draw (6.473,12.015) circle (1pt);
\draw (6.551,2.032) circle (1pt);
\draw (6.551,11.968) circle (1pt);
\draw (6.631,11.920) circle (1pt);
\draw (6.631,2.080) circle (1pt);
\draw (6.712,2.130) circle (1pt);
\draw (6.712,11.870) circle (1pt);
\draw (6.792,11.819) circle (1pt);
\draw (6.792,2.181) circle (1pt);
\draw (6.873,11.767) circle (1pt);
\draw (6.873,2.233) circle (1pt);
\draw (6.952,11.713) circle (1pt);
\draw (6.952,2.287) circle (1pt);
\draw (7.000,7.000) circle (1pt);
\draw (7.032,2.342) circle (1pt);
\draw (7.032,11.658) circle (1pt);
\draw (7.111,11.601) circle (1pt);
\draw (7.111,2.399) circle (1pt);
\draw (7.189,11.543) circle (1pt);
\draw (7.189,2.457) circle (1pt);
\draw (7.267,2.517) circle (1pt);
\draw (7.267,11.483) circle (1pt);
\draw (7.344,11.421) circle (1pt);
\draw (7.344,2.579) circle (1pt);
\draw (7.420,2.642) circle (1pt);
\draw (7.420,11.358) circle (1pt);
\draw (7.496,2.707) circle (1pt);
\draw (7.496,11.293) circle (1pt);
\draw (7.505,7.000) circle (1pt);
\draw (7.570,11.226) circle (1pt);
\draw (7.570,2.774) circle (1pt);
\draw (7.644,11.156) circle (1pt);
\draw (7.644,2.844) circle (1pt);
\draw (7.716,2.915) circle (1pt);
\draw (7.716,11.085) circle (1pt);
\draw (7.786,2.988) circle (1pt);
\draw (7.786,11.012) circle (1pt);
\draw (7.855,10.937) circle (1pt);
\draw (7.855,3.063) circle (1pt);
\draw (7.871,7.000) circle (1pt);
\draw (7.922,3.139) circle (1pt);
\draw (7.922,10.861) circle (1pt);
\draw (7.959,7.000) circle (1pt);
\draw (7.959,7.012) circle (1pt);
\draw (7.959,6.988) circle (1pt);
\draw (7.960,7.025) circle (1pt);
\draw (7.960,6.975) circle (1pt);
\draw (7.961,6.963) circle (1pt);
\draw (7.961,7.037) circle (1pt);
\draw (7.962,6.951) circle (1pt);
\draw (7.962,7.049) circle (1pt);
\draw (7.964,7.062) circle (1pt);
\draw (7.964,6.938) circle (1pt);
\draw (7.966,7.074) circle (1pt);
\draw (7.966,6.926) circle (1pt);
\draw (7.968,6.914) circle (1pt);
\draw (7.968,7.086) circle (1pt);
\draw (7.971,6.901) circle (1pt);
\draw (7.971,7.099) circle (1pt);
\draw (7.974,7.111) circle (1pt);
\draw (7.974,6.889) circle (1pt);
\draw (7.977,7.123) circle (1pt);
\draw (7.977,6.877) circle (1pt);
\draw (7.981,6.864) circle (1pt);
\draw (7.981,7.136) circle (1pt);
\draw (7.985,6.852) circle (1pt);
\draw (7.985,7.148) circle (1pt);
\draw (7.988,3.218) circle (1pt);
\draw (7.988,10.782) circle (1pt);
\draw (7.990,7.161) circle (1pt);
\draw (7.990,6.839) circle (1pt);
\draw (7.995,7.174) circle (1pt);
\draw (7.995,6.826) circle (1pt);
\draw (8.001,7.186) circle (1pt);
\draw (8.001,6.814) circle (1pt);
\draw (8.007,6.801) circle (1pt);
\draw (8.007,7.199) circle (1pt);
\draw (8.013,6.788) circle (1pt);
\draw (8.013,7.212) circle (1pt);
\draw (8.020,7.225) circle (1pt);
\draw (8.020,6.775) circle (1pt);
\draw (8.028,7.239) circle (1pt);
\draw (8.028,6.761) circle (1pt);
\draw (8.036,6.748) circle (1pt);
\draw (8.036,7.252) circle (1pt);
\draw (8.045,6.734) circle (1pt);
\draw (8.045,7.266) circle (1pt);
\draw (8.051,3.298) circle (1pt);
\draw (8.051,10.702) circle (1pt);
\draw (8.054,7.279) circle (1pt);
\draw (8.054,6.721) circle (1pt);
\draw (8.065,7.293) circle (1pt);
\draw (8.065,6.707) circle (1pt);
\draw (8.076,7.307) circle (1pt);
\draw (8.076,6.693) circle (1pt);
\draw (8.088,7.321) circle (1pt);
\draw (8.088,6.679) circle (1pt);
\draw (8.101,6.664) circle (1pt);
\draw (8.101,7.336) circle (1pt);
\draw (8.112,3.379) circle (1pt);
\draw (8.112,10.621) circle (1pt);
\draw (8.115,6.649) circle (1pt);
\draw (8.115,7.351) circle (1pt);
\draw (8.130,7.366) circle (1pt);
\draw (8.130,6.634) circle (1pt);
\draw (8.146,7.381) circle (1pt);
\draw (8.146,6.619) circle (1pt);
\draw (8.164,7.396) circle (1pt);
\draw (8.164,6.604) circle (1pt);
\draw (8.172,10.538) circle (1pt);
\draw (8.172,3.462) circle (1pt);
\draw (8.184,7.412) circle (1pt);
\draw (8.184,6.588) circle (1pt);
\draw (8.205,6.572) circle (1pt);
\draw (8.205,7.428) circle (1pt);
\draw (8.228,3.547) circle (1pt);
\draw (8.228,10.453) circle (1pt);
\draw (8.229,6.556) circle (1pt);
\draw (8.229,7.444) circle (1pt);
\draw (8.249,7.000) circle (1pt);
\draw (8.256,6.539) circle (1pt);
\draw (8.256,7.461) circle (1pt);
\draw (8.263,8.393) circle (1pt);
\draw (8.263,5.607) circle (1pt);
\draw (8.275,5.604) circle (1pt);
\draw (8.275,8.396) circle (1pt);
\draw (8.282,3.633) circle (1pt);
\draw (8.282,10.367) circle (1pt);
\draw (8.285,6.523) circle (1pt);
\draw (8.285,7.477) circle (1pt);
\draw (8.295,8.403) circle (1pt);
\draw (8.295,5.597) circle (1pt);
\draw (8.318,6.506) circle (1pt);
\draw (8.318,7.494) circle (1pt);
\draw (8.325,8.413) circle (1pt);
\draw (8.325,5.587) circle (1pt);
\draw (8.333,10.281) circle (1pt);
\draw (8.333,3.719) circle (1pt);
\draw (8.356,6.490) circle (1pt);
\draw (8.356,7.510) circle (1pt);
\draw (8.365,8.430) circle (1pt);
\draw (8.365,5.570) circle (1pt);
\draw (8.381,3.807) circle (1pt);
\draw (8.381,10.193) circle (1pt);
\draw (8.400,6.474) circle (1pt);
\draw (8.400,7.526) circle (1pt);
\draw (8.417,5.541) circle (1pt);
\draw (8.417,8.459) circle (1pt);
\draw (8.425,3.895) circle (1pt);
\draw (8.425,10.105) circle (1pt);
\draw (8.451,7.540) circle (1pt);
\draw (8.451,6.460) circle (1pt);
\draw (8.464,3.985) circle (1pt);
\draw (8.464,10.015) circle (1pt);
\draw (8.477,5.496) circle (1pt);
\draw (8.477,8.504) circle (1pt);
\draw (8.500,9.924) circle (1pt);
\draw (8.500,4.076) circle (1pt);
\draw (8.511,7.549) circle (1pt);
\draw (8.511,6.451) circle (1pt);
\draw (8.531,9.831) circle (1pt);
\draw (8.531,4.169) circle (1pt);
\draw (8.536,5.429) circle (1pt);
\draw (8.536,8.571) circle (1pt);
\draw (8.559,4.263) circle (1pt);
\draw (8.559,9.737) circle (1pt);
\draw (8.583,6.451) circle (1pt);
\draw (8.583,7.549) circle (1pt);
\draw (8.583,5.945) circle (1pt);
\draw (8.583,8.055) circle (1pt);
\draw (8.584,4.358) circle (1pt);
\draw (8.584,9.642) circle (1pt);
\draw (8.586,8.655) circle (1pt);
\draw (8.586,5.345) circle (1pt);
\draw (8.606,9.547) circle (1pt);
\draw (8.606,4.453) circle (1pt);
\draw (8.622,5.249) circle (1pt);
\draw (8.622,8.751) circle (1pt);
\draw (8.625,4.549) circle (1pt);
\draw (8.625,9.451) circle (1pt);
\draw (8.642,4.646) circle (1pt);
\draw (8.642,9.354) circle (1pt);
\draw (8.646,8.851) circle (1pt);
\draw (8.646,5.149) circle (1pt);
\draw (8.655,9.256) circle (1pt);
\draw (8.655,4.744) circle (1pt);
\draw (8.660,5.047) circle (1pt);
\draw (8.660,8.953) circle (1pt);
\draw (8.662,7.533) circle (1pt);
\draw (8.662,6.467) circle (1pt);
\draw (8.663,9.156) circle (1pt);
\draw (8.663,4.844) circle (1pt);
\draw (8.665,4.945) circle (1pt);
\draw (8.665,9.055) circle (1pt);
\draw (8.742,6.504) circle (1pt);
\draw (8.742,7.496) circle (1pt);
\draw (8.812,6.560) circle (1pt);
\draw (8.812,7.440) circle (1pt);
\draw (8.869,7.373) circle (1pt);
\draw (8.869,6.627) circle (1pt);
\draw (8.912,7.298) circle (1pt);
\draw (8.912,6.702) circle (1pt);
\draw (8.942,6.779) circle (1pt);
\draw (8.942,7.221) circle (1pt);
\draw (8.959,7.000) circle (1pt);
\draw (8.959,6.858) circle (1pt);
\draw (8.959,7.142) circle (1pt);
\draw (8.963,6.935) circle (1pt);
\draw (8.963,7.065) circle (1pt);
\end{tikzpicture}
\end{center}
\caption{Complex chromatic roots of $X_{H,\,W_4}(100)$}
\label{fig:complex1}
\end{figure}
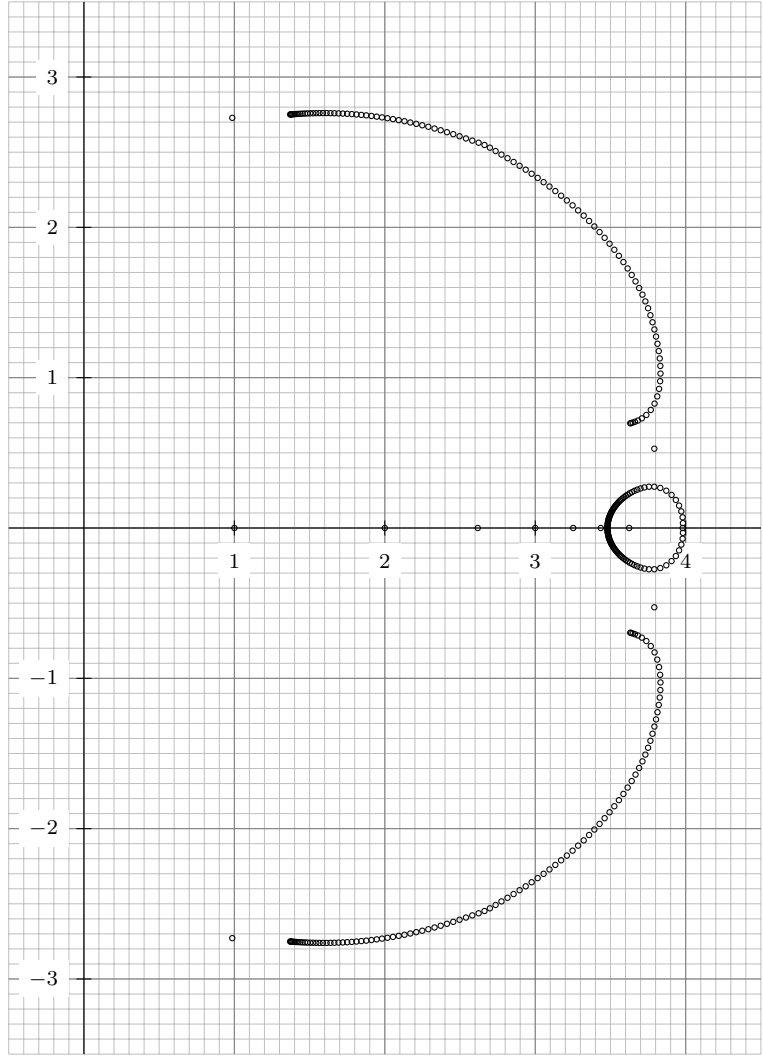

For these polynomials, it is also possible to compute their {\em complex} chromatic roots; the roots for $X_{H,\,W_4}(100)$ are plotted in Figure~\ref{fig:complex1} and are in close agreement with the theoretical limiting curves as first shown in Beraha \& Kahane \cite{MR539072}.

The computation of the complex chromatic roots starts to become time consuming for larger chromatic polynomials, but we can go much further if we are only interested a single real root. By repeated squaring, we can easily compute $(MD)^2$, $(MD)^4$, $\ldots$, $(MD)^{512}$ and the bisection method (with rational arithmetic) can be used to locate a real chromatic root near 4, which is presumably the largest. The values found in this manner are given in Table~\ref{tab:reals2}. 

\begin{table}
\begin{center}
\begin{tabular}{|c|c|}
\hline
$n$&Value\\
\hline
2&3.848343257\\
4&3.875604098\\
8& 3.905148525\\
16& 3.932717391\\
32& 3.955237394\\
64&3.971732100\\
128&3.982848013\\
256&3.989898687\\
512&3.994181944\\
\hline
\end{tabular}
\end{center}
\caption{Real chromatic roots of $X_{H,\,W_4}(n+1)$ close to 4.}
\label{tab:reals2}
\end{table}

\subsection{Notes on the computations}

As the results of this section are computational in nature, it is important to verify them as much as feasible.  To this end, the partitioned chromatic polynomial of the base graph $H$ was computed by two independent programs, being Gary Haggard's fast {\tt chromial}  program \cite{MR1675926} and a slower but much simpler program of the author. The computation of the chromatic polynomials was performed symbolically by Maple and then using exact arithmetic we confirmed that the polynomial satisfied Tutte's identity for planar triangulations, i.e that 
$$
P(G,\tau+2) = (\tau+2) \tau^{3n-10} P(G,\tau+1)^2
$$
where $\tau$ is the golden ratio $(1+\sqrt5)/2$ and $n$ is the  number of vertices of 
$G$. The computation of the complex chromatic roots of these polynomials was done with Bini and Fiorentino's {\tt MPSolve}  package \cite{MR1772050} which rigorously controls the output accuracy by dynamically altering the precision during computation. The computation of the largest real root was then cross-checked by using the bisection method with exact rational arithmetic in Maple. 

\section{Roots Tending to Four}

To prove that a family of graphs $X_{A,B}(n)$ has real chromatic roots arbitrarily close to four, we will show that for any sufficiently small $\eps > 0$, there is an integer $n_0$ such that for all $n > n_0$
$$
P(X_{A,B}(n), 4-\eps) < 0 \quad {\rm and} \quad P(X_{A,B}(n), 4) > 0.
$$
We are mainly concerned with the situation when $X_{A,B}(n)$ is planar, in which case the second inequality automatically holds.

The result of taking a matrix to a high power is essentially determined by its spectral properties. In this case, if $0 < \eps < 2-\tau$, then for $x=4-\eps$, the matrix
$MD$ has four distinct positive eigenvalues
\begin{align*}
\lambda_1&=2\\
\lambda_2&=2 - 5 \eps + 10 \eps^2 / 3 + O(\eps^3)\\
\lambda_3&=2 - 8 \eps + 26 \eps^2 / 3 + O(\eps^3)\\
\lambda_4 &=0
\end{align*}
with corresponding (right) eigenvectors $v_1 = (1,-1,-1,1)^T$, 
{
 $$v_2=
\left(\begin{array}{r}3/2+ 35\eps/12 + 1103\eps^2 / 216+O(\eps^3) \\ 
1 + 5\eps/3 + 85 \eps^2/27 + O(\eps^3) \\
1 + 5\eps/3 + 85 \eps^2/27 + O(\eps^3) \\
\hfill -1\hfill\end{array}\right),
$$
$$
v_3 = 
\left(\begin{array}{r}\phantom{XX }- \eps/3 - 14\eps^2 / 27+O(\eps^3) \\ 
1/2 + \eps/6 + 4 \eps^2/27 + O(\eps^3) \\
1/2 + \eps/6 + 4 \eps^2/27 + O(\eps^3) \\
\hfil1\hfil\end{array}\right)
$$
}and $v_4 = (0,1,-1,0)^T$. The matrix $MD$ is not symmetric, and so the eigenvectors are not necessarily orthogonal. However, they {\em do} form an orthogonal set as long as we are willing to use a different inner product, namely
$$
\langle v, w\rangle = v^T D w,
$$
where here, and subsequently, we identify the $1 \times 1$ matrix with the scalar it contains.

\begin{lemma}\label{lem:orth}
If $\{v_1,\ldots,v_4\}$ are the eigenvectors of $MD$ for some fixed $x=4-\eps$ where 
$\eps \in (0,2-\tau)$ then 
$$
v_i^T D v_j = 0 \iff i \not= j.
$$
\end{lemma}

\proof
The entries of $D$ are all strictly positive and $v_i \not= 0$ so it follows that $v_i^T D v_i$ is a positive linear combination of squares with at least one non-zero term, and so 
$v_i^T D v_i > 0$. Now as $v_i$ and $v_j$ are eigenvectors of $MD$ we have 
\begin{align*}
(v_i^T D) (MD v_j)&= \lambda_j v_i^T D v_j,\\
(v_j^T D) (MD v_i)&= \lambda_i v_j^T D v_i.
\end{align*}
The second expression is simply the transpose of the first, and so
$$
\lambda_j v_i^T D v_j = \lambda_i v_i^T D v_j.
$$
If $i\not= j$, then $\lambda_i \not= \lambda_j$, and so it must be the case that
$v_i^T D v_j = 0$.
\qed

Let $\|v\|$ denote the norm of the vector $v$ with respect to this inner product. As the eigenvalues $\{\lambda_i\}$ are distinct, the eigenvectors $\{v_i\}$ are linearly independent, so we can express $Q(A)$ and $Q(B)$ as linear combinations 
\begin{align*}
Q(A) &= \alpha_1 v_1 + \alpha_2 v_2 + \alpha_3 v_3 + \alpha_4 v_4,\\
Q(B) &= \beta_1 v_1 + \beta_2 v_2 + \beta_3 v_3 + \beta_4 v_4.
\end{align*} 
Then the value of 
$Q(A)^T D (MD)^{n-1} Q(B)$ at $x=4-\eps$ is simply
$$
\alpha_1 \beta_1 \lambda_1^{n-1} \|v_1\|^2  + \alpha_2\beta_2 \lambda_2^{n-1} \|v_2\|^2 + \alpha_3\beta_3 \lambda_3^{n-1} \|v_3\|^2 + \alpha_4\beta_4 \lambda_4^{n-1}  \|v_4\|^2.
$$

As $n$ increases, this sum is dominated by the largest eigenvalue $\lambda_i$ for which {\em both} $\alpha_i$ and $\beta_i$ are non-zero, and the eventual {\em sign} of this sum depends purely on whether the dominant $\alpha_i$ and $\beta_i$ have the same or opposite signs.

We can now prove the result implied by the numerical computations of the previous section.

\begin{theorem}
Graphs in the family $X_{H,\,W_4}(n)$ have real chromatic roots 
arbitrarily close to 4.
\end{theorem}
\proof
As above, suppose that $x=4-\eps$ with $\eps \in (0,2-\tau)$, that 
$\{v_1,\ldots,v_4\}$ are the eigenvectors for $MD$, and that 
$Q(H) = \sum \alpha_i v_i$ and $Q(W_4) = \sum \beta_i v_i$. Then
substituting $x=4-\eps$ into $Q(H)^T D v_i$ and $Q(W_4)^TD v_i$ for $i=1,\,2$ and expressing the results as series in $\eps$, we get
\begin{align*}
\alpha_1 \|v_1\|^2&= 0,\\
\beta_1 \|v_1\|^2&= 0,\\
\alpha_2 \|v_2 \|^2&= -50 \eps + O(\eps^2),\\
\beta_2 \|v_2\|^2 &= 5 + 20 \eps / 3 + O(\eps^2).\\
\end{align*}
Therefore the dominant eigenvalue for this family is $\lambda_2$ and for sufficiently small $\eps>0$, the values $\alpha_2$ and $\beta_2$ have opposite signs. \qed

\section{Other end-graphs}

It is natural to consider which other end-graphs can be used in this construction to produce real chromatic roots tending to 4. The graphs $X_{A,B}(n)$ are planar if and only if both $A$ and $B$ are planar with a face of size 4 as the distinguished 4-cycle. In this situation, the coefficient of $v_1$ in the expressions for $Q(A)$ and $Q(B)$ is always zero.

\begin{lemma}
If $A$ is a planar graph with a face $a_1a_2a_3a_4$ of size 4 as distinguished 4-cycle, then 
$$
Q(A)^T D v_1 = 0.
$$
\end{lemma}

\proof
As $v_1 = (1,-1,-1,1)^T$ we need to prove that 
$$
P_1(A,x)/\ff2 + P_4(A,x)/\ff4  = P_2(A,x)/\ff3 + P_3(A,x)/\ff3.
$$
If $x$ is an integer, then each term of the form $P_i(A,x)/\ff s$ is the number of $x$-colourings of $A$ of type $i$ divided by the number of ways of choosing the $s$ colours used on $\{a_1,\ldots,a_4\}$. Equivalently, each term is the number of ways of extending a particular colouring of $\{a_1,\ldots,a_4\}$ to an $x$-colouring of the entire graph. 
Divide the colourings of $A$ using colours from $1$, $2$, $\ldots\,$, $x$ into two classes as follows: $\varphi$ is in Class 1 if  
$$[\varphi(a_1), \varphi(a_2), \varphi(a_3), \varphi(a_4)] = [1,2,1,2]\ {\rm or}\  [1,2,3,4],$$
and Class 2 if
$$[\varphi(a_1), \varphi(a_2), \varphi(a_3), \varphi(a_4)] = [1,2,1,3]\ {\rm or}\ [1,2,3,2].$$

We will now find a bijection between the two classes, and thus prove that the stated expression is true for all integers, and therefore for all $x$.

Let $\varphi$ be any colouring in either class, and consider whether or not there is a $2$--$4$ path (i.e. a path whose vertices are all coloured 2 or 4) connecting $a_2$ and $a_4$. If there is a $2$--$4$ path between $a_2$ and $a_4$, then there is {\em not} a $1$--$3$ path connecting $a_1$ and $a_3$ because $a_1a_2a_3a_4$ is a face of a planar graph. Define a new colouring $\varphi'$ by applying the following rules: if there is a $2$--$4$ path from $a_2$ to $a_4$, then exchange the colours $1$ and $3$ on the $1$--$3$ component of $A$ containing $a_3$, and if there is no $2$--$4$ path from $a_2$ to $a_4$ then exchange colours $2$ and $4$ on the $2$--$4$ component of $A$ containing $a_4$. It is straightforward to check that $\varphi'$ is in Class 1 if and only $\varphi$ is in Class 2, and vice versa. \qed

We note that if $A$ is planar, but the distinguished 4-cycle is not a face, then the conclusion of the above lemma may or may not hold. 

If $C$ is a planar graph with a distinguished face of size 4, then say that it is {\em positive} if the expression for $Q(C)^TDv_2$ at $x=4-\eps$ has a positive leading term when expressed as a series in $\eps$, and {\em negative} otherwise. Then the family $X_{A,B}(n)$ has real chromatic roots arbitrarily close to 4 if one of $A$ and $B$ is positive, and the other negative.

\begin{lemma}\label{lem:alpha1}
If $A$ is a planar graph with a face $a_1a_2a_3a_4$ of size 4 as distinguished 4-cycle, then the constant term in the expression for $\alpha_2 \|v_2\|^2$ is $5/24$ times the number of 4-colourings of $A / a_1a_3 / a_2a_4$ (i.e the graph obtained from $A$ by identifying $a_1$ with $a_3$ and $a_2$ with $a_4$).
\end{lemma}
\proof
When $\eps = 0$, the eigenvector $v_2 = (3/2, 1, 1, -1)^T = (5/2, 0,0,0)^T + (-1,1,1,-1)^T$ and so $Q(A)^T D v_2 = 5/2\ P_1(A,4) / 12 = 5/24\ P_1(A,4)$.
\qed

\begin{figure}
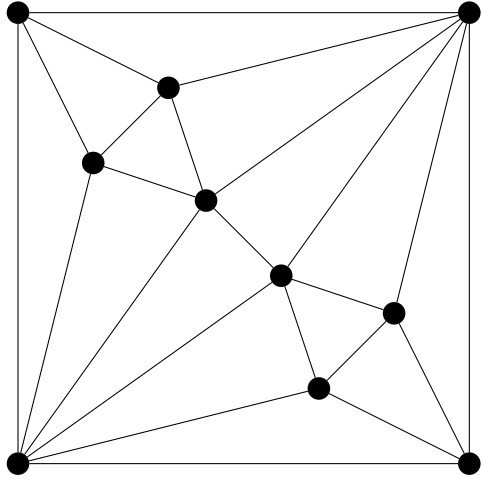

\begin{center}
\begin{pgfpicture}{0cm}{0cm}{6cm}{6cm}
\pgfnodecircle{v0}[fill]{\pgfxy(0,0)}{0.15cm}
\pgfnodecircle{v1}[fill]{\pgfxy(6,0)}{0.15cm}
\pgfnodecircle{v2}[fill]{\pgfxy(6,6)}{0.15cm}
\pgfnodecircle{v3}[fill]{\pgfxy(0,6)}{0.15cm}
\pgfnodeconnline{v0}{v1}
\pgfnodeconnline{v1}{v2}
\pgfnodeconnline{v2}{v3}
\pgfnodeconnline{v3}{v0}
\pgfnodecircle{v4}[fill]{\pgfxy(4,1)}{0.15cm}
\pgfnodecircle{v5}[fill]{\pgfxy(5,2)}{0.15cm}
\pgfnodecircle{v6}[fill]{\pgfxy(2,5)}{0.15cm}
\pgfnodecircle{v7}[fill]{\pgfxy(1,4)}{0.15cm}
\pgfnodeconnline{v0}{v4}
\pgfnodeconnline{v4}{v1}
\pgfnodeconnline{v1}{v5}
\pgfnodeconnline{v5}{v2}
\pgfnodeconnline{v2}{v6}
\pgfnodeconnline{v6}{v3}
\pgfnodeconnline{v3}{v7}
\pgfnodeconnline{v7}{v0}
\pgfnodecircle{v8}[fill]{\pgfxy(3.5,2.5)}{0.15cm}
\pgfnodecircle{v9}[fill]{\pgfxy(2.5,3.5)}{0.15cm}
\pgfnodeconnline{v4}{v5}
\pgfnodeconnline{v6}{v7}
\pgfnodeconnline{v8}{v9}
\pgfnodeconnline{v4}{v8}
\pgfnodeconnline{v0}{v8}
\pgfnodeconnline{v0}{v9}
\pgfnodeconnline{v5}{v8}
\pgfnodeconnline{v2}{v8}
\pgfnodeconnline{v2}{v9}
\pgfnodeconnline{v6}{v9}
\pgfnodeconnline{v7}{v9}
\end{pgfpicture}
\end{center}
\caption{A 10-vertex negative planar graph}
\label{fig:small}
\end{figure}

From this lemma, it follows that a {\em necessary} (but not sufficient) condition for a planar graph with distinguished 4-face to be negative is that the graph obtained by contracting the two diagonals of the face has no 4-colourings. Searching among the triangulations of a square with this property, negative graphs with as few as 10 vertices were found, with one example shown in Figure~\ref{fig:small}.

\subsection{Notes}

When Jacobsen, Salas and Sokal \cite{MR2000228} examined the cylindrical triangular lattice without end-graphs, they obtained a $2 \times 2$ transfer matrix, with only $\lambda_2$ and $\lambda_3$ as eigenvalues. The reason for this is that rather than using colour partitions to index the rows and columns of the transfer matrix, they exploited planarity by using only ``non-crossing non-nearest neighbour'' partitions. 
This is essentially equivalent to avoiding consideration of partitions of Type 1 which, as we have just shown, contain no additional information when $A$ is planar. They gained another dimension by exploiting rotational symmetry, thus working directly in the subspace orthogonal to $v_4$. Therefore the eigenvalues $\lambda_1=2$ and $\lambda_4=0$ do not occur in their analysis.

Ro{\v{c}}ek, Shrock \& Tsai \cite{roshts} obtained the same four eigenvalues for the transfer matrix as above, and also noticed that for certain pairs of end-graphs, the 
dominant eigenvalue is $\lambda_1 = 2$ (or rather, $\lambda_4=2$ in their notation)
and that this caused the limiting curves of chromatic roots to change. In particular, the small closed curve passing through $z=4$ and $z \simeq 3.481$ becomes much larger, and passes through $z=4$ and $z=1+\tau \simeq 2.618$ instead (Fig. 4(b) of \cite{roshts}).

\section{Acknowledgements}

Thanks are due to Douglas Woodall for finding the original record-holding graph and sending it to me, Gary Haggard for his {\tt chromial} program that allowed me to convince myself that this graph really was very special, Alan Sokal and Jes\'us Salas for a number of ideas about how to generalize this single data point, and subsequent helpful discussions on this topic, and Robert Shrock for pointers to the physics literature on this topic. Graphics were created with Till Tantau's Ti{\em k}Z and PGF packages \cite{tikz}.

\bibliographystyle{plain}
\bibliography{plantri}

\end{document}